\documentclass[11pt]{article}
\usepackage{amsmath,amssymb}
\usepackage{fullpage}
\usepackage[dvips]{graphicx}
\usepackage{psfrag}  
\def\COMMENT#1{}
\def\TASK#1{}
\def\noproof{{\unskip\nobreak\hfill\penalty50\hskip2em\hbox{}\nobreak\hfill%
        $\square$\parfillskip=0pt\finalhyphendemerits=0\par}\goodbreak}
\def\endproof{\noproof\bigskip}
\newdimen\margin   
\def\textno#1&#2\par{%
    \margin=\hsize
    \advance\margin by -4\parindent
           \setbox1=\hbox{\sl#1}%
    \ifdim\wd1 < \margin
       $$\box1\eqno#2$$%
    \else
       \bigbreak
       \hbox to \hsize{\indent$\vcenter{\advance\hsize by -3\parindent
       \sl\noindent#1}\hfil#2$}%
       \bigbreak
    \fi}
\def\proof{\removelastskip\penalty55\medskip\noindent{\bf Proof. }}
\newtheorem{firstthm}{Proposition}
\newtheorem{thm}[firstthm]{Theorem}
\newtheorem{prop}[firstthm]{Proposition}
\newtheorem{fact}[firstthm]{Fact}
\newtheorem{lemma}[firstthm]{Lemma}

\newtheorem{conj}[firstthm]{Conjecture}
\newtheorem{claim}[firstthm]{Claim}
\newtheorem{ques}[firstthm]{Question} 

\title{On perfect packings in dense graphs}
\author{J\'ozsef Balogh,\footnote{University of Illinois, Urbana-Champaign, USA and University of California, San Diego, USA, jobal@math.uiuc.edu.
This author is supported by NSF CAREER Grant DMS-0745185, UIUC Campus Research
Board Grant 11067, and OTKA Grant K76099.}
\ Alexandr V. Kostochka\footnote{
 University of Illinois,
Urbana-Champaign, USA and Institute of Mathematics, Novosibirsk,
Russia, kostochk@math.uiuc.edu. This author
is supported in part by NSF grant DMS-0965587 and by grant 12-01-00631
of the Russian Foundation for Basic Research.}
\  and Andrew Treglown\footnote{Queen Mary, University of London, United Kingdom, treglown@maths.qmul.ac.uk}}
\begin{document}
\label{firstpage} \maketitle

\begin{abstract} We say that a graph $G$ has a perfect $H$-packing if there exists a set of vertex-disjoint 
copies of $H$ which cover all the vertices in $G$. We consider various problems concerning perfect $H$-packings:
Given $n, r , D \in \mathbb N$, we characterise the edge density threshold that ensures a perfect $K_r$-packing in any graph
$G$ on $n$ vertices and with minimum degree $\delta (G) \geq D$. We also give two conjectures concerning degree sequence conditions
which force a graph to contain a perfect $H$-packing. Other related embedding problems are also considered.
Indeed, we give a structural result concerning $K_r$-free graphs that satisfy
a certain degree sequence condition.
\end{abstract}

\section{Introduction} 
Given two graphs $H$ and $G$, a \emph{perfect $H$-packing} in $G$ is a collection of vertex-disjoint copies of $H$ which cover
all the vertices in $G$. Perfect $H$-packings are also referred to as \emph{$H$-factors} or \emph{perfect $H$-tilings}. 
Hell and Kirkpatrick~\cite{hell} showed that the decision problem
whether a graph $G$ has a perfect $H$-packing is NP-complete precisely when $H$ has a
component consisting of at least $3$ vertices. So for such graphs $H$, it is unlikely that there is a complete characterisation
of those graphs containing a perfect $H$-packing. Thus, there has been significant attention on obtaining sufficient conditions that
ensure a graph $G$ contains a perfect $H$-packing.

A seminal result in the area is the Hajnal-Szemer\'edi theorem~\cite{hs}
which states that a graph $G$ whose order $n$
is divisible by $r$ has a perfect $K_r$-packing provided that $\delta (G) \geq (r-1)n/r$. 
K\"uhn and Osthus~\cite{kuhn, kuhn2} characterised, up to an additive constant, the minimum degree which ensures a graph $G$ 
contains a perfect $H$-packing for an arbitrary graph $H$. 

It is easy to see that the minimum degree condition in the Hajnal-Szemer\'edi theorem cannot be lowered. Of course, this does
not mean that one cannot strengthen this result.
\emph{Ore-type} degree conditions consider the sum of the degrees of non-adjacent vertices in a graph.
The following Ore-type result of Kierstead and Kostochka~\cite{kier} implies the Hajnal-Szemer\'edi theorem.
\begin{thm}[Kierstead and Kostochka~\cite{kier}]\label{orekk}
Let $n,r \in \mathbb N$ such that $r  $ divides $n$. Suppose that $G$ is a graph on $n$ vertices such that for all non-adjacent
$x \not =y \in V(G)$,
$$d(x)+d(y) \geq 2(1-1/r)n-1 .$$
Then $G$ contains a perfect $K_r$-packing.
\end{thm}
K\"uhn, Osthus and Treglown~\cite{kotore}  characterised, asymptotically, the Ore-type degree condition which ensures a graph $G$ 
contains a perfect $H$-packing for an arbitrary graph $H$.

\subsection{Perfect packings in dense graphs of low minimum degree}
It is easy to characterise the edge density that forces a graph $G$ to contain a perfect $K_r$-packing when there are no other
restrictions.
Indeed, given $n,r \in \mathbb N$ such that $r \geq 2$ divides $n$, if $G$ is a graph on $n$ vertices and
$e(G) \geq \binom{n}{2}- n+r$ then $G$ contains a perfect $K_r$-packing. Moreover, if $G$ is a copy $K$ of $K_{n-1}$ together with a
vertex which sends precisely $r-2$ edges to $K$, then $e(G)=\binom{n}{2}- n+r-1$ and $G$ does not contain a perfect $K_r$-packing.
The following result of Akiyama and Frankl~\cite{frank} refines this observation.
\begin{thm}[Akiyama and Frankl~\cite{frank}]\label{frank}
Let $n,r \in \mathbb N$ such that $r$ divides $n$. Suppose $G$ is a graph on $n$ vertices and
$e(\overline{G}) \leq \min \{ \binom{n/r+1}{2} , n-r+1 \}$. Then $G$ has a perfect $K_r$-packing
unless $\overline{G}$ is isomorphic to one of the following graphs:
\begin{itemize}
\item[(i)] A copy of $K_{n/r+1}$ together with $(1-1/r)n-1$ isolated vertices;
\item[(ii)] The disjoint union of $K_{1,n-r-j+1}$, $j$ edges and $r-j-2$ isolated vertices, for some $1\leq j \leq r-2$.
\end{itemize}
\end{thm}
When (for example) $n \geq r^3$, $\binom{n/r+1}{2} > n-r+1$. Hence, in this case Theorem~\ref{frank} is 
equivalent to the following: If $G$ is a graph on $n$ vertices and $e(G)\geq \binom{n}{2}- n+r-1$ then either 
$G$ contains a perfect $K_r$-packing or $\overline{G}$ is isomorphic to a graph as in (ii).

In Sections~\ref{edgeexsec} and~\ref{edgesec} we consider the following natural problem: Let $n,r \in \mathbb N$ such that $r$ divides $n$. Given some
$D \in \mathbb N$, what edge density condition ensures that any graph $G$ on $n$ vertices and 
of minimum degree $\delta (G) \geq D$ contains
a perfect $K_r$-packing? 

We fully resolve the problem, and our answers for $r=2$ and $r\geq 3$ differ.

\begin{thm}\label{p1}
For an even positive $n$ and integer $1\leq d<n/2$, let 
$h(n,d):=\binom{n-d-1}{2}+d(d+1)$ and let
$f(2,n,d)$ denote the maximum integer $c$ such that  some $n$-vertex graph
with minimum degree at least $d$ and at least $c$ edges has no perfect matching. Then 
\begin{equation*}
f(2,n,d)=\max\{h(n,d),h(n,0.5n-1)\}.
\end{equation*}\end{thm}

\COMMENT{Is $r=2$ worth thinking about/trivial?}
\begin{thm}\label{mainthm1}
Let $n,r \in \mathbb N$ such that $r\geq 3$ and $r$ divides $n$. Given any $D \in \mathbb N$ such that
$r-1\leq D\leq (r-1)n/r-1$ define
$$g(n,r,D):= \max \left \{ \binom{n}{2}-\binom{n/r+1}{2}, D(n-D)+\binom{n-1-D}{2}+e({T}(D, r-2))     \right \}.$$
Suppose that $G$ is a graph on $n$ vertices with $\delta (G) \geq D$ and $e(G)>g(n,r,D)$. Then $G$ contains a perfect $K_r$-packing.
Moreover, there exists a graph $G'$ on $n$ vertices with $\delta (G')\geq D$ and $e(G')=g(n,r,D)$ but such that $G'$ does not
contain a perfect $K_r$-packing.
\end{thm}
Clearly a graph $G$ of minimum degree $\delta (G) <r-1$ cannot contain a perfect $K_r$-packing. Further,
regardless of edge density, every graph $G$ whose order $n$ is divisible by $r$ and with $\delta (G)\geq (r-1)n/r$
contains a perfect $K_r$-packing. Thus, Theorem~\ref{mainthm1} covers all values of $D$ where our problem was not solved
previously.

\medskip

An \emph{equitable $k$-colouring} of a graph $G$ is a proper $k$-colouring of $G$ such that any two colour classes differ in size by
at most one. Let $n,r  \in \mathbb N$ such that $r$ divides $n$. Notice that a graph $G$ on $n$ vertices has a perfect $K_r$-packing
if and only if the complement $\overline{G}$ of $G$ has an equitable $n/r$-colouring. So, for example, the Hajnal-Szemer\'edi
theorem can be stated in terms of equitable colourings: Let $G$ be a graph on $n$ vertices such that $r$ divides $n$. If $\Delta(G)
\leq n/r-1$ then $G$ has an equitable $n/r$-colouring.

It is often easier to work in the language of equitable colourings compared to perfect packings. Indeed, rather than
prove Theorem~\ref{orekk} directly,
Kierstead and Kostochka proved the equivalent statement for equitable colourings. Here we also find
it more convenient to work with equitable colourings. Thus, instead of proving Theorem~\ref{mainthm1} directly we prove the
following equivalent result.

\begin{thm}\label{t1}
Let $n,r \in \mathbb N$ such that $r\geq 3$ and $r$ divides $n$. Recall that $T(n,r)$ denotes the Tur\' an graph.
Given any $D \in \mathbb N$ such that
$n/r\leq D\leq n-r$ define
$$f(n,r,D):= \min \left \{  \binom{n/r+1}{2} , D+e(\overline{T}(n-D-1, r-2))     \right \}.$$
Suppose that $G$ is a graph on $n$ vertices with $\Delta (G) \leq D$ and $e(G)<f(n,r,D)$. Then $G$ has an equitable $n/r$-colouring.
Moreover, there exists a graph $G'$ on $n$ vertices with $\Delta (G')\leq D$ and $e(G')=f(n,r,D)$ but such that $G'$ does not
have an equitable $n/r$-colouring.
\end{thm}

We prove Theorem~\ref{p1} and describe extremal constructions for Theorems~\ref{mainthm1} and~\ref{t1} in Section~\ref{edgeexsec}.
That is, we show that the edge density condition in Theorem~\ref{mainthm1} is best possible for all  values of $D$
such that $r-1\leq D \leq (r-1)n/r-1$.
Section~\ref{edgesec} contains a proof of Theorem~\ref{t1}.

\subsection{Degree sequence conditions forcing a perfect packing}
Chv\'atal~\cite{ch} gave a condition on the degree sequence of a graph which ensures Hamiltonicity:
Suppose that $G$ is a graph on $n$ vertices and that the degrees of the graph are $d_1 \leq \dots \leq d_n$. 
If $n \geq 3$ and $d_i \geq i+1$ or $d_{n-i} \geq n-i$ for all $i<n/2$ then G is Hamiltonian.
The following is a simple consequence of Chv\'atal's theorem.
\begin{thm}[Chv\'atal~\cite{ch}]\label{hampath}
Suppose that $G$ is a graph on $n\geq 2$ vertices and the degrees of the graph are $d_1\leq \dots \leq d_n$.
If $$d_i \geq i \ \text{ or } \ d_{n-i+1}\geq n-i \ \text{ for all } \ 1 \leq i \leq n/2$$
then $G$ contains a Hamilton path.
\end{thm}

We propose the following conjecture on the degree sequence of a graph which forces a perfect $K_r$-packing.
\begin{conj}\label{conj1}
Let $n,r \in \mathbb N$ such that $r $ divides $n$. Suppose that $G$ is a graph on $n$ vertices with degree sequence $d_1 \leq 
\dots \leq d_n$ such that:
\begin{itemize}
\item[($\alpha$)] $d_i \geq (r-2)n/r+i $ for all $i < n/r$;
\item[($\beta$)] $d_{n/r+1} \geq (r-1)n/r$.
\end{itemize}
Then $G$ contains a perfect $K_r$-packing.
\end{conj}
\COMMENT{Note that this is slightly stronger than my original conjecture.}
Note that Conjecture~\ref{conj1}, if true, is much stronger than the Hajnal-Szemer\'edi theorem since the degree condition allows
for $n/r$ vertices to have degree less than $(r-1)n/r$. 
Further, Proposition~\ref{extremal1} in Section~\ref{secextremal}
 shows that the condition on the degree sequence in Conjecture~\ref{conj1} 
is essentially ``best possible''.
It is easy to see that Theorem~\ref{hampath} implies Conjecture~\ref{conj1} in the case when $r=2$. \COMMENT{Chvatal had a typo in the statement of his result. But Theorem~\ref{hampath} is what he meant!}
We prove the conjecture in the case when $G$ is additionally $K_{r+1}$-free
(see Section~\ref{sec2}).

If one can prove Conjecture~\ref{conj1}, it seems likely it can be used to prove the next conjecture.
\begin{conj}\label{conj2} Suppose $\gamma >0$ and $H$ is a graph with $\chi (H)=r$. Then there exists an integer $n_0=n_0 (\gamma ,H)$
such that the following holds. If $G$ is a graph whose order $n \geq n_0$ is divisible by $|H|$, and whose degree
sequence $d_1\leq \dots \leq d_n$ satisfies
\begin{itemize}
\item  $d_i \geq (r-2)n/r+i +\gamma n  $ for all $i < n/r$,
\end{itemize}
then $G$ contains a perfect $H$-packing. 
\end{conj}
Since first submitting this paper, the third author and~Knox~\cite{knox} have proven Conjecture~\ref{conj2} in the case when
$r=2$. (In fact, they have proven a much more general result concerning embedding spanning bipartite graphs of small bandwidth.)

The following result of Erd\H{o}s~\cite{erdos1} characterises those degree sequences which force a copy of $K_r$ in a graph $G$.
\begin{thm}[Erd\H{o}s~\cite{erdos1}]
Let $G$ be a graph on $n$ vertices with degree sequence $d_1\leq \dots \leq d_n$. If $G$ is $K_{r+1}$-free then there is an $r$-partite
graph $G'$ on $n$ vertices whose degree sequence $d'_1\leq \dots \leq d'_n$ satisfies
$$d_i \leq d'_i \ \text{ for all } \ i \leq n.$$
\end{thm}
In Section~\ref{turansec} we prove the following related structural theorem.
\begin{thm}\label{turan}
Suppose that $n,r \in \mathbb N$ such that $n \geq r$ and so that $r$ divides $n$.
Let $G$ be a $K_{r+1}$-free graph on $n$ vertices whose degree sequence $d_1 \leq \dots \leq d_n$
is such that
$d_{ n/r } \geq (r-1)n/r$.
Then $G \subseteq T(n,r)$, where  $T(n,r)$ is 
the complete $r$-partite \emph{Tur\'an graph} on $n$ vertices; so each vertex class has size $\lceil n/r 
\rceil$ or $\lfloor n/r \rfloor$.
\end{thm}

\section{The case $r=2$ and extremal examples for $r\geq 3$}\label{edgeexsec}

\subsection{Perfect matchings in dense graphs}\label{sectutte}
In this section we establish the density threshold that ensures every graph $G$ on an even number $n$ of vertices and of minimum
degree $\delta (G) \geq d$ contains a perfect matching. Note that we only consider values of $d$ such that $1 \leq d <n/2$, since if 
$\delta (G)\geq n/2 $ then $G$ has a perfect matching, regardless of the edge density.

Recall that 
$h(n,d):=\binom{n-d-1}{2}+d(d+1).$
 Note that for a fixed even $n$, $h(n,d)$ decreases with $d$ in the interval
$[0,n/3-5/6]$ and increases with $d$ in $[n/3-5/6,0.5n-1]$.

For a positive even $n$ and an integer $0\leq d<n/2$, let
$A,B$ and $C$ be disjoint sets with  $|A|=d+1$, $|B|=d$, $|C|=n-2d-1$.
Let $H=H(n,d)$ be the graph with the vertex set
$A\cup B\cup C$ such that $H[B\cup C]=K_{n-d-1}$, and each vertex in $A$
is adjacent to each vertex in $B$ and to no vertex in $C$.
So $H$ does not contain a perfect matching and has exactly $h(n,d)$ edges.

 The examples of $H(n,d)$ show that 
$f(2,n,d)\geq \max\{h(n,d),h(n,0.5n-1)\}$. 
Thus to derive Theorem~\ref{p1}, it  suffices to prove that   an $n$-vertex graph $G$ with $ \delta(G) \geq d$ and 
$e(G) >\max\{h(n,d),$ $h(n,0.5n-1)\}$ has a perfect matching.

Consider such a graph $G$. Let $d_1 \leq \dots \leq d_n$ denote the degree sequence of $G$.
If $d_i \geq i$ for all $1 \leq i \leq n/2$ then Theorem~\ref{hampath} implies that $G$ contains a perfect matching.
Suppose for a contradiction that $d_{i'} \leq i'-1$ for some $1 \leq i' \leq n/2$. Note that $i' >d$ as $\delta (G) \geq d$.

Let $A$ denote the set of $i'$ vertices in $G$ that correspond to the first $i'$ terms $d_1 , \dots , d_{i'}$ of the degree sequence.
Set $B:=V(G) \backslash A$. Then
$$e(G[B]) \geq e(G) - i'(i'-1) > \max\{h(n,d),h(n,0.5n-1)\} - i'(i'-1)$$ since $d(x) \leq i'-1$ for all $x \in A$.
Note that $\max\{h(n,d),h(n,0.5n-1)\} \geq h(n,i'-1)$ since $d<i' \leq n/2$.
Therefore, 
$$e(G[B])>\max\{h(n,d),h(n,0.5n-1)\} - i'(i'-1) \geq h(n,i'-1) -i'(i'-1){=} \binom{n-i'}{2},$$
a contradiction as $|B|=n-i'$. Thus, $d_i \geq i$ for all $1 \leq i \leq n/2$, as desired.

\subsection{Examples for $r\geq 3$}

We will give the extremal examples for Theorem~\ref{t1}. Since Theorems~\ref{mainthm1} and~\ref{t1} are equivalent,
the complements of the extremal graphs for Theorem~\ref{t1} are the extremal graphs for Theorem~\ref{mainthm1}.
\begin{prop}\label{propex1} Suppose that $n,r \in \mathbb N$ such that $r\geq 3$ and $r$ divides $n$.  
Then there exists a graph $G_1$ on $n$ vertices such that $\Delta (G_1)=n/r$,
$$e(G_1)=\binom{n/r+1}{2},$$
but such that $G_1$ does not have an equitable $n/r$-colouring.
\end{prop}
\proof
Let $G_1$ denote the disjoint union of a clique $V$ on $n/r+1$ vertices and an independent set $W$ of $(1-1/r)n-1$ vertices. 
So every independent set in $G_1$ contains at most one vertex from $V$. But since $|V|=n/r+1$,
$G_1$ does not have an  equitable $n/r$-colouring. Further, $\Delta (G_1)=n/r$ and
$e(G_1)=\binom{n/r+1}{2}.$
\endproof

\begin{prop}\label{propex2} Suppose that $n,r \in \mathbb N$ such that $r\geq 3$ and $n=kr$ for some $k \geq 2$. Further,
let $D \in \mathbb N$ such that $n/(r-1) \leq D \leq n-r$.  
Then there exists a graph $G_2$ on $n$ vertices such that $\Delta (G_2)=D$,
$$e(G_2)=D+e(\overline{T}(n-D-1, r-2)) ,$$
but such that $G_2$ does not have an equitable $n/r$-colouring.
\end{prop}
\proof Let $G_2$ denote the disjoint union of a copy $K$ of $K_{1,D}$ and a copy of $\overline{T}(n-D-1, r-2)$. So $|G|=n$.
Let $v$ denote the vertex of degree $D$ in $K$. The largest independent set in $G_2$ that contains $v$ is of size $r-1$. Thus,
$G_2$ does not have an equitable $n/r$-colouring. Further, $e(G_2)=D+e(\overline{T}(n-D-1, r-2)).$ 

Since $n/(r-1)\leq D$ we have that $n-1 \leq (r-1)D$. Thus, every vertex in the copy of 
$\overline{T}(n-D-1, r-2)$ has degree at most
$$\left \lceil \frac{n-D-1}{r-2} \right \rceil -1 \leq \frac{n-D-1}{r-2} \leq D.$$ This implies that $\Delta (G_2)=D$.
\endproof
\COMMENT{Awkward that the example in Prop~\ref{propex2} does not have $\Delta (G_2)=D$ if $D<n/(r-1)$. This is why we need the
following extra proposition.}
Clearly Propositions~\ref{propex1} and~\ref{propex2} show that one cannot lower the edge density condition in Theorem~\ref{t1}
in the case when $n/(r-1) \leq D \leq n-r$. The following result, together with Proposition~\ref{propex1}, shows that Theorem~\ref{t1} is best possible in the case
when $n/r \leq D \leq n/(r-1)$.
\begin{prop}\label{propex3} Let $n,r \in \mathbb N$ such that $r\geq 3$ and $r$ divides $n\geq 2r$. 
Suppose that $D\in \mathbb N$ such that $n/r\leq D \leq n/(r-1)$. Then
$$f(n,r,D)=\binom{n/r+1}{2}.$$
\end{prop}
The following simple consequence of Tur\'an's theorem will be used in the proof of Theorem~\ref{t1}.
\begin{fact}\label{turanbound}
Let $n,r\in \mathbb N$ such that $r\leq n$. Then $$e(T(n,r))\leq \left(1-\frac{1}{r}\right)\frac{n^2}{2} \ \text{ and thus }  \
e(\overline{T}(n,r))\geq \frac{n^2}{2r}-\frac{n}{2}.$$
\end{fact}
We will also require the following easy result.
\begin{lemma}\label{case3prop}
Let $n,r \in \mathbb N$ such that $r\geq 4$ and $r$ divides $n \geq 3r$. Suppose that $D \in \mathbb N$ such that
$n/r \leq D < (n+r)/(r-1)$. Then 
$$f(n,r,D)=\binom{n/r+1}{2}.$$
\end{lemma}

\section{Proof of Theorem~\ref{t1}}\label{edgesec}

\subsection{Preliminaries}
Suppose for a contradiction that the result is false. Let $G$ be a counterexample with the fewest vertices. That is, $n=|V(G)|=rk$
for some $k \in \mathbb N$, $\Delta (G) \leq D$ for some $D \in \mathbb N$ such that $n/r\leq D\leq n-r$, $e(G)<f(n,r,D)$ and
$G$ has no equitable $n/r$-colouring. By the Hajnal-Szemer\'edi theorem, $\Delta (G)\geq n/r$.
Notice that given fixed $n$ and $r$, $f(n,r,D)$ is non-increasing with respect to $D$. Thus, we may assume that
$\Delta (G)=D$.

We first show that $k \geq 4$. Indeed, if $n=2r$ then $f(n,r,D)\leq\binom{3}{2}=3$. But it is easy to see that
every graph $G_1$ on $2r$ vertices and with $e(G_1)\leq 2$ has an equitable $2$-colouring. If $n=3r$ then
$f(n,r,D)\leq \binom{4}{2}=6$. Consider any graph $G_1$ on $3r$ vertices with $e(G_1)\leq 5$ and $3\leq \Delta (G_1)\leq 5$.
Let $x$ denote the vertex in $G_1$ where $d_{G_1} (x)=\Delta (G_1)$. Since $3\leq d_{G_1}(x)\leq 5$, $x$ lies 
in an independent set $I$ in $G_1$ of size $r$. But then $G_1-I$ contains $2r$ vertices and at most $2$ edges. So $G_1-I$ has an
equitable $2$-colouring and hence $G_1$ has an equitable $3$-colouring. 

Let $v \in V(G)$ such that $d_G (v)=D$. Set $G^*:=G-(N_G (v) \cup \{v\})$.
Since $f(n,r,D)\leq D+e(\overline{T}(n-D-1, r-2)) $ we have that $e(G^*)<e(\overline{T}(n-D-1, r-2))$. Thus, by Tur\'an's theorem, $G^*$ 
contains an independent set of size $r-1$. Hence, $v$ lies in an independent set in $G$ of size $r$. 
Amongst all such independent sets of size $r$ that contain $v$, choose a set $I=\{v, x_1, \dots , x_{r-1} \}$ such that
$d_G (x_1)+\dots + d_G (x_{r-1})$ is maximised.

Set $G':=G-I$, $n':=|V(G')|=n-r$ and $D':=\Delta (G')\leq D$. Notice that $D'\geq n'/r$. 
(Indeed, if not, then by the Hajnal-Szemer\'edi
theorem $G'$ contains an equitable $n'/r$-colouring. Thus, as $I$ is an independent set in $G$ this gives us an equitable 
$n/r$-colouring of $G$, a contradiction.)
Furthermore, $D' \leq n'-r$. If not then 
\begin{align*}
e(G)\geq D+D' \geq 2D' \geq 2(n'-r+1)=2n-4r+2
\end{align*}
and further,
\begin{align*}
e(G)& <f(n,r,D)\leq f(n,r, n-2r+1) \leq (n-2r+1) +e(\overline{T} (2r-2,r-2)) \\
& \leq (n-2r+1)+(r+3)=n-r+4.
\end{align*} 
\medskip  Therefore, $2n-4r+2 <n-r+4$ and so $n<3r+2$ a contradiction since $n =kr \geq  4r$.

Since $n'/r\leq D' \leq n'-r$, if  $e(G')< f(n',r,D')$ then the minimality of $G$ implies that $G'$ has an equitable
$n'/r$-colouring. This then implies that $G$ has an equitable $n/r$-colouring, a contradiction. Thus,
\begin{align}\label{contra1}
e(G')\geq f(n',r,D').
\end{align}
We now split our argument into three cases.

\subsection{Case 1: $ \bf f(n',r, D')=\binom{n'/r+1}{2}$. }
By~(\ref{contra1}), $e(G')\geq \binom{n'/r+1}{2}=\binom{n/r}{2}$. Since $d_G (v)=D\geq n/r$,
$$e(G)\geq \frac{n}{r}+\binom{n/r}{2}=\binom{n/r+1}{2} \geq f(n,r,D),$$
a contradiction, as desired.

\subsection{Case 2: $\bf D'\leq D-1$ and $ \bf f(n',r, D')=D'+ e(\overline{T}(n'-D'-1,r-2))$. }
The following claim will be useful.
\begin{claim}\label{claim1}
$D'<\frac{r-1}{2r-3} n -\frac{(r^2-r+1)}{2r-3}.$
\end{claim}
\proof
Note that 
\begin{align}\label{claimeq1}
D+D'+e(\overline{T}(n'-D'-1,r-2))\stackrel{(\ref{contra1})}{\leq} e(G)<f(n,r,D)\leq D+e(\overline{T}(n-D-1,r-2)).
\end{align}
Since $D'\leq D-1$, clearly $e(\overline{T}(n'-D,r-2)) \leq e(\overline{T}(n'-D'-1,r-2))$. Thus, (\ref{claimeq1}) implies
that 
\begin{align}\label{claimeq2}
D'+e(\overline{T}(n'-D,r-2)) < e(\overline{T}(n-D-1,r-2)).
\end{align}
One can obtain $\overline{T}(n-D-1,r-2)$ from $\overline{T}(n'-D,r-2)$ by adding $r-1$ vertices and at most
\begin{align}\label{claimeq3}
(n'-D)+\frac{n-D-2}{r-2} \ \text{ edges.}
\end{align}
Hence~(\ref{claimeq2})~and~(\ref{claimeq3}) give
$$D'<n'-D+\frac{n-D-2}{r-2}.$$
Rearranging, and using that $D'\leq D-1$ and $n'=n-r$ we get that
$$\left(2+\frac{1}{r-2}\right)D'< \left(1+\frac{1}{r-2}\right)n- \frac{(r^2-r+1)}{r-2}.$$
Thus,
$$D'<\frac{r-1}{2r-3} n -\frac{(r^2-r+1)}{2r-3},$$
as desired.
\endproof
Since we are in Case~2 we have that
\begin{align}\label{EQ1}
D'+e(\overline{T}(n-r-D'-1,r-2))\leq \binom{n'/r+1}{2}=\binom{n/r}{2}.
\end{align}
Notice that for fixed $n$ and $r$, $D'+e(\overline{T}(n-r-D'-1,r-2))$ is non-increasing as $D'$ increases. Hence,
(\ref{EQ1}) and Claim~\ref{claim1} imply that 
\begin{align}\label{EQ2}
D''+e(\overline{T}(n-r-D''-1,r-2))\leq \frac{n^2}{2r^2}-\frac{n}{2r}
\end{align}
where $D'':=\lfloor (r-1)n/(2r-3)  -(r^2-r+1)/(2r-3) \rfloor$.
Note that $$n-r-\frac{r-1}{2r-3}n  +\frac{(r^2-r+1)}{2r-3}-1=\frac{r-2}{2r-3}n+\frac{4-r^2}{2r-3}.$$
So Fact~\ref{turanbound} and (\ref{EQ2}) imply that
\begin{align*}
& \left(\frac{r-1}{2r-3}n-  \frac{(r^2-r+1)}{2r-3} -\frac{(2r-4)}{2r-3} \right)    
+\frac{1}{2(r-2)}\left( \frac{r-2}{2r-3} n+\frac{4-r^2}{2r-3}\right) ^2\\ & -\frac{1}{2} 
\left( \frac{r-2}{2r-3} n+\frac{4-r^2}{2r-3}\right) \leq \frac{n^2}{2r^2}-\frac{n}{2r} .
\end{align*}

Next we will move all terms from the previous equation to the left hand side and simplify. The coefficient of $n^2$ is
\begin{align}\label{n2}
\frac{r-2}{2(2r-3)^2}-\frac{1}{2r^2}=\frac{r^3-6r^2+12r-9}{2r^2(2r-3)^2}.
\end{align}
The coefficient of $n$ is
\begin{align}\label{n1}
\frac{r-1}{2r-3}-\frac{(r-2)}{2(2r-3)}+\frac{1}{2r}+\frac{(4-r^2)}{(2r-3)^2}=\frac{r^2-4r+9}{2r(2r-3)^2}.
\end{align}
The constant term is
\begin{align}\label{c}
-\frac{(r^2+r-3)}{2r-3}+\frac{(r^2-4)^2}{2(r-2)(2r-3)^2}+\frac{(r^2-4)}{2(2r-3)}=\frac{-r^4+3r^3+4r^2-26r+28}{2(r-2)(2r-3)^2}.
\end{align}
\COMMENT{Coefficients of $n^2$ and $n$ are always non-negative. So we can indeed substitute $n\geq 4r$ into (\ref{eqcon1}).}
Since $n \geq 4r$, (\ref{n2})--(\ref{c}) imply that
\begin{align}\label{eqcon1}
\frac{8(r^3-6r^2+12r-9)}{(2r-3)^2}+\frac{2(r^2-4r+9)}{(2r-3)^2}+\frac{-r^4+3r^3+4r^2-26r+28}{2(r-2)(2r-3)^2} \leq 0.
\end{align}
Multiplying (\ref{eqcon1}) by $2(r-2)(2r-3)^2$ we get
\begin{align*}
15r^4 -121 r^3 +364 r^2 -486r +244 \leq 0
\end{align*}
This yields a contradiction, since it is easy to check that
\begin{align*}
15r^4 -121 r^3 +364 r^2 -486r +244 > 0
\end{align*}
for all $r \in \mathbb N$ such that $r \geq 3$.

\subsection{Case 3: $\bf D'=D$ and $ \bf f(n',r, D')=D'+ e(\overline{T}(n'-D'-1,r-2))$.}
By (\ref{contra1}) we have that
\begin{align}\label{thresh1}
e(G') \geq f(n',r,D')=D'+ e(\overline{T}(n'-D'-1,r-2)).
\end{align}
Consider any vertex $x \in V(G')$ such that $d_{G'} (x)=D'=D$. Since $\Delta (G)=D$, $x$ is not adjacent to any vertex 
in $I=\{v, x_1, \dots , x_{r-1} \}$. Further, $I$ was chosen such that $d_G (x_1)+\dots +d_G (x_{r-1})$ is maximised. Thus,
 $d_G(x_1)=\dots =d_G (x_{r-1})=D$.
Together with (\ref{thresh1}) this implies that
\begin{align}\label{thresh2}
e(G)\geq (r+1)D +e(\overline{T}(n'-D-1,r-2)).
\end{align}
Since $e(G)<f(n,r,D)\leq D+ e(\overline{T} (n-D-1,r-2))$, (\ref{thresh2}) implies that
\begin{align}\label{thresh3}
rD+e(\overline{T}(n'-D-1,r-2)) <e(\overline{T}(n-D-1,r-2)).
\end{align}
One can obtain $\overline{T}(n-D-1,r-2)$ from $\overline{T}(n'-D-1,r-2)$ by adding $r$ vertices and at most
\begin{align}\label{thresh4}
(n'-D-1)+\frac{2(n-D-3)}{r-2}+1 \ \text{ edges.}
\end{align}
Thus, (\ref{thresh3}) and (\ref{thresh4}) imply that
$$rD<n-r-D + \frac{2(n-D-3)}{r-2}$$ 
and so 
\begin{align}\label{thresh5}
\left(r+1+\frac{2}{r-2}\right)D< \left(1+\frac{2}{r-2}\right)n +\frac{(-r^2+2r-6)}{r-2}<\left( 1+\frac{2}{r-2} \right)n.
\end{align}
If $r=3$ then (\ref{thresh5}) implies that
$$D<\frac{n}{2}.$$
Since $f(n',3,D)=\min\{\binom{n'/3+1}{2}, D+\binom{n'-D-1}{2} \}$
it is easy to see that if $f(n',3,D)=D+\binom{n'-D-1}{2}$ then $D\geq 2n'/3+1=2n/3-1$.
Thus, $2n/3-1 \leq D <n/2$, a contradiction since $n \geq 4r=12$.

If $r\geq 4$  then (\ref{thresh5}) implies that
$$D<\frac{n}{r-1}=\frac{n'}{r-1}+\frac{r}{r-1}. $$
Since $n'\geq 3r$, Lemma~\ref{case3prop} implies that $f(n',r,D')=\binom{n'/r+1}{2}$ and so we are in Case~1, which we have already
dealt with.

\section{The extremal examples for Conjecture~\ref{conj1}}\label{secextremal}
\begin{prop}\label{extremal1}
Suppose that $n,r,k \in \mathbb N$ such that $r\geq 2$ divides $n$ and $1\leq k <n/r$. Then there exists a graph $G$ on $n$ vertices whose
degree sequence $d_1\leq \dots \leq d_n$ satisfies
\begin{itemize}
\item $d_i = (r-2)n/r+k-1$ for all $1 \leq i \leq k$;
\item $d_i = (r-1)n/r$ for all $k+1 \leq i \leq (r-2)n/r+k$;
\item $d_i = n-k-1$ for all $(r-2)n/r+k+1 \leq i \leq n-k+1$;
\item $d_i =n-1$ for all $n-k+2 \leq i \leq n$,
\end{itemize}
but such that $G$ does not contain a perfect $K_r$-packing.
\end{prop}
\proof Let $G'$ denote the complete $(r-2)$-partite graph whose vertex classes $V_1, \dots , V_{r-2}$ each have size $n/r$. Obtain $G$ from $G'$ by adding
the following vertices and edges:
Add a set $V_{r-1}$ of $2n/r-2k+1$ vertices to $G'$, a set $V_r$ of $k-1$ vertices and a set $V_0$ of $k$ vertices. Add all edges from $V_0 \cup V_{r-1} \cup V_r$ to $V_1 \cup \dots \cup V_{r-2}$. Further, add all edges with both endpoints in $V_{r-1} \cup V_{r}$. Add all possible edges between $V_0$ and $V_r$. 

So $V_0$ is an independent set, and there are no edges between $V_0$ and $V_{r-1}$. This implies that any copy of $K_r$ in $G$ containing
a vertex from $V_0$ must also contain at least one vertex from $V_r$. But since $|V_0|>|V_{r}|$ this implies that $G$ does not
contain a perfect $K_r$-packing. Furthermore, $G$ has our desired degree sequence.
\endproof
\COMMENT{Add picture!}
Notice that the graphs $G$ considered in Proposition~\ref{extremal1} satisfy ($\beta$) from Conjecture~\ref{conj1} and only
fail to satisfy ($\alpha$) in the case when $i=k$ (and in this case $d_k =(r-2)n/r+k-1$). 
\medskip

Let $n,r \in \mathbb N$ such that $r$ divides $n$. Denote by $T^*(n,r)$ the complete $r$-partite graph on $n$ vertices 
 with $r-2$ vertex classes of
size $n/r$, one vertex class of size $n/r-1$ and one vertex class of size $n/r+1$. Then $T^*(n,r)$ does not
contain a perfect $K_r$-packing. 
Furthermore, $T^*(n,r)$ satisfies ($\alpha$) but condition ($\beta$) fails; we have that 
$d _{n/r+1}=(r-1)n/r-1$ here.
Thus, together $T^*(n,r)$ and Proposition~\ref{extremal1} show that, if true, Conjecture~\ref{conj1} is essentially best possible.

\section{A special case of Conjecture~\ref{conj1}}\label{sec2}
We now give a simple proof of Conjecture~\ref{conj1} in the case when $G$ is $K_{r+1}$-free.
\begin{thm}\label{Krfree}
Let $n,r \in \mathbb N$ such that $r \geq2 $ divides $n$. Suppose that $G$ is a graph on $n$ vertices with degree sequence $d_1 \leq \dots \leq d_n$ such that:
\begin{itemize}
\item $d_i \geq (r-2)n/r+i $ for all $i < n/r$;
\item $d_{n/r+1} \geq (r-1)n/r$.
\end{itemize}
Further suppose that no vertex $x \in V(G)$ of degree less than $(r-1)n/r$ lies in a copy of $K_{r+1}$.
Then $G$ contains a perfect $K_r$-packing.
\end{thm}
\proof
We prove the theorem by induction on $n$. In the case when $n=r$ then $d_{n/r+1}=d_2\geq (r-1)r/r=r-1$. 
This implies that every vertex in $G$ has degree $r-1$.
Hence $G=K_r$ as desired.
So suppose that $n>r$ and the result holds for smaller values of $n$. Let $x_1 \in V(G)$ such that $d_G(x_1)=d_1 \geq (r-2)n/r+1$.
If $d_G (x_1) \geq (r-1)n/r$ then $\delta (G)\geq (r-1)n/r$. Thus $G$ contains a perfect $K_r$-packing by the Hajnal-Szemer\'edi
theorem. So we may assume that $(r-2)n/r+1\leq d_G (x_1) <(r-1)n/r$. In particular, $x_1$ does not lie in a copy of $K_{r+1}$.
We first find a copy of $K_r$ containing $x_1$. If $r=2$, $x_1$ has a neighbour and so we have our desired copy of $K_2$. So assume
that $r\geq 3$.
Certainly $N_G(x_1)$ contains a vertex $x_2$ such that $d_G(x_2) \geq (r-1)n/r$. Thus, $|N_G(x_1)\cap N_G(x_2)|\geq (r-3)n/r+1 >0$. 
So if $r=3$ we obtain our desired copy of $K_r$. Otherwise, we can find a vertex $x_3 \in N_G(x_1)\cap N_G(x_2)$ such that $d_G(x_3) \geq (r-1)n/r$. We can repeat this argument until we have obtained vertices $x_1, \dots , x_r$ that together form a copy $K'_r$ of $K_r$.

Let $G':=G-V(K'_r)$ and set $n':=n-r=|V(G')|$. 
Since $G$ does not contain a copy of $K_{r+1}$ containing $x_1$, every vertex $x \in V(G)\backslash V(K'_r)$
sends at most $r-1$ edges to $K'_r$ in $G$.  Thus, $d_{G'} (x) \geq d_G (x) -(r-1)$ for 
all $x \in V(G')$. So if $d_G (x) \geq (r-1)n/r$ then $d_{G'} (x) \geq (r-1)n/r -(r-1)=(r-1)n'/r$ for all $x \in V(G')$. If a vertex
$y \in V(G')$ does not lie in a copy of $K_{r+1}$ in $G$ then clearly $y$ does not lie in a copy of $K_{r+1}$ in $G'$. This means 
that no vertex $y \in V(G')$ of degree less than $(r-1)n'/r$ lies in a copy of $K_{r+1}$.

Let $d'_1 \leq \dots \leq d'_{n'}$ denote the degree sequence of $G'$. It is easy to check that
$d'_i \geq {(r-2)n'}/{r}+i$ for all $i < n'/r$ and that $d'_{n'/r+1} \geq (r-1)n'/r$. 
Indeed, since $x_1 \in V(K'_r)$ where $d_G (x_1)=d_1$, we have that $d'_i \geq d_{i+1}-(r-1)$ for all $1 \leq i \leq n'$.
Thus, for all $1 \leq i < n'/r=n/r-1$, $d'_i \geq d_{i+1}-(r-1) \geq (r-2)n/r+(i+1)-(r-1)=(r-2)n'/r+i$. Similarly,
$d'_{n'/r+1}=d'_{n/r} \geq d_{n/r+1}-(r-1)\geq (r-1)n/r-(r-1)=(r-1)n'/r$.
Hence, by induction $G'$ contains a perfect $K_r$-packing. Together with 
$K'_r$ this gives us our desired perfect $K_r$-packing in $G$.
\endproof

\section{Proof of Theorem~\ref{turan}}\label{turansec}
Consider any $x_1 \in V(G)$ such that $d_G(x_1)\geq (r-1)n/r$.
Since $d_{n/r} \geq (r-1)n/r$ we can greedily select vertices $x_2, \dots , x_{r-1}$ such that
\begin{itemize}
\item $x_1, \dots , x_{r-1}$ induce a copy of $K_{r-1}$ in $G$;
\item $d_G(x_i) \geq (r-1)n/r$ for all $1 \leq i \leq r-1$.
\end{itemize}
Note that since $G$ is $K_{r+1}$-free, $\cap ^{r-1} _{i=1} N_G (x_i)$ is an independent set. The choice of $x_1, \dots ,$ $x_{r-1}$ 
implies that
$|\cap ^{r-1} _{i=1} N_G (x_i)|\geq n/r$. Let $V_1$ denote a subset of $\cap ^{r-1} _{i=1} N_G (x_i)$ of size $n/r$. Thus
$V_1$ contains a vertex $x^1 _1$ of degree at least $(r-1)n/r$. 

As before we can find vertices $x^1_2, \dots , x^1_{r-1}$ such that
\begin{itemize}
\item $x^1 _1, \dots , x^1 _{r-1}$ induce a copy of $K_{r-1}$ in $G$;
\item $d_G(x^1_i) \geq (r-1)n/r$ for all $1 \leq i \leq r-1$.
\end{itemize}
So $\cap ^{r-1} _{i=1} N_G (x^1 _i)$ is an independent set of size at least $n/r$. 
Let $V_2$ denote a subset of $\cap ^{r-1} _{i=1} N_G (x^1 _i)$ of size $n/r$. Note that $N_G (x^1 _1) \cap V_1=\emptyset$
since $x^1 _1 \in V_1$
and $V_1$ is an independent set. Thus as $V_2 \subseteq N_G (x^1 _1)$, $V_1 \cap V_2 =\emptyset$.

Our aim is to find disjoint sets $V_1, \dots , V_r \subseteq V(G)$ of size $n/r$ and vertices
$x^1 _1, \dots, x^1 _{r-1},$ $x^2 _1, \dots, x^2 _{r-1}$, $\dots , x^{r-1} _1,
\dots, x^{r-1} _{r-1}$ with the following properties:
\begin{itemize}
\item $G[V_j]$ is an independent set for all $1 \leq j \leq r$;
\item Given any $1\leq j \leq r-1$, $x^j _k \in V_k$ for each $1 \leq k \leq j$;
\item $d_G (x^j _k) \geq (r-1)n/r$ for all $1 \leq j \leq r-1$ and $1 \leq k \leq r-1$;
\item $x^j _1, \dots , x^j _{r-1}$ induce a copy of $K_{r-1}$ in $G$ for all $1 \leq j \leq r-1$.
\end{itemize}
Clearly finding such a partition $V_1, \dots , V_r$ of $V(G)$ implies that $G \subseteq T(n,r)$.

Suppose that for some $1< j <r$ we have defined sets $V_1, \dots , V_j$ and vertices
\\$x^1 _1, \dots, x^1 _{r-1}, \dots , x^{j-1} _1, \dots, x^{j-1} _{r-1}$ with our desired properties.
Since $d_{n/r} \geq (r-1)n/r$ and $V_1, \dots , V_j$ are independent sets of size $n/r$ we can choose vertices
$x^j _1, \dots, x^j _j$ such that for all $1 \leq k \leq j$,
\begin{itemize}
\item $x^j _k \in V_k$ and $d_G (x^j _k) \geq (r-1)n/r$.
\end{itemize}
This degree condition, together with the fact that $x^j _1, \dots, x^j _j$ lie in different vertex classes, implies that
these vertices form a copy of $K_j$ in $G$. We now greedily select further vertices $x^j_{j+1}, \dots , x^j_{r-1}$ such that
\begin{itemize}
\item $x^j _1, \dots , x^j _{r-1}$ induce a copy of $K_{r-1}$ in $G$;
\item $d_G(x^j_k) \geq (r-1)n/r$ for all $j+1 \leq k \leq r-1$.
\end{itemize}
So $\cap ^{r-1} _{i=1} N_G (x^j _i)$ is an independent set of size at least $n/r$. 
Let $V_{j+1}$ denote a subset of $\cap ^{r-1} _{i=1} N_G (x^j _i)$ of size $n/r$. Note that, for each $1 \leq k \leq j$,
 $N_G (x^j _k) \cap V_{k}=\emptyset$
since $x^j _k \in V_k$
and $V_k$ is an independent set. Thus as $V_{j+1} \subseteq N_G (x^j _k)$ for each $1 \leq k \leq j$, $V_{j+1}$ is disjoint from 
$V_1 \cup \dots \cup V_j$.
 
Repeating this argument we obtain our desired
sets $V_1, \dots , V_r \subseteq V(G)$ and vertices
$x^1 _1, \dots, x^1 _{r-1},$ $x^2 _1, \dots, x^2 _{r-1} , \dots , x^{r-1} _1, \dots, x^{r-1} _{r-1}$.

\section{Possible extensions of Conjecture~\ref{conj1}}
We suspect that the following `Chv\'atal-type' degree sequence condition forces a graph to contain a perfect $K_r$-packing.
\begin{ques}\label{ques1}
Let $n,r \in \mathbb N$ such that $r\geq 2 $ divides $n$. Suppose that $G$ is a graph on $n$ vertices with degree sequence $d_1 \leq 
\dots \leq d_n$ such that for all $i \leq  n/r$:
\begin{itemize}
\item $d_i \geq (r-2)n/r+i $ or $d_{n-i(r-1)+1} \geq n-i$.
\end{itemize}
Does this condition imply that $G$ contains a perfect $K_r$-packing?
\end{ques}
\COMMENT{Ques slightly different to original one (originally $d_{n-i(r-1)}$). But this seems the correct condition.}
Note that Theorem~\ref{hampath} answers this question in the affirmative when $r=2$.
The following example shows that we cannot have a lower value in the second part of the condition in Question~\ref{ques1}. 
\begin{prop}\label{extremal2}
Suppose that $n,r,k \in \mathbb N$ such that $r\geq 2$ divides $n$ and $1\leq k \leq n/r$. Then there exists a graph $G$ on $n$ vertices whose
degree sequence $d_1\leq \dots \leq d_n$ satisfies
\begin{itemize}
\item $d_{n-i(r-1)+1} \geq n-i$ for all $i \in [n/r]\backslash \{k\}$;
\item $d_{n-k(r-1)+1}=n-k-1$,
\end{itemize} 
but such that $G$ does not contain a perfect $K_r$-packing.
\end{prop}
\proof Let $G$ be the graph on $n$ vertices with vertex classes $V_1, V_2$ and $V_3$ of sizes $k$, $(r-1)k-1$ and $n-rk+1$ 
respectively and with the following edges: There are all possible edges between $V_1$ and $V_2$ and between $V_2$ and $V_3$.
Further add all possible edges in $V_2$ and all edges in $V_3$. Thus, $V_1$ is an independent set and there are no edges between
$V_1$ and $V_3$.

The degree sequence of $G$ is
$$\underbrace{(r-1)k-1, \dots ,(r-1)k-1}_{k \text{ times}}, \underbrace{n-k-1, \dots , n-k-1}_{n-rk+1
\text{ times}}, \underbrace{n-1, \dots , n-1}_{(r-1)k-1 \text{ times}}.$$
Hence $G$ satisfies our desired degree sequence condition. Every copy $K'_r$ of $K_r$ in $G$ that contains a vertex from $V_1$ must
contain $r-1$ vertices from $V_2$. But since $|V_1|(r-1)>|V_2|$ this implies that $G$ does not contain a perfect $K_r$-packing.
\endproof
\COMMENT{add picture}

The $r$th power of a Hamilton cycle
$C$ is obtained from $C$ by adding an edge between every pair of vertices of distance at most $r$ on $C$. Seymour~\cite{seymour} conjectured the
following strengthening of Dirac's theorem.
\begin{conj}[P\'osa-Seymour, see~\cite{seymour}]\label{seymour}
Let $G$ be a graph on $n$ vertices. If $\delta (G) \geq \frac{r}{r+1} n$
then $G$ contains the $r$th power of a Hamilton cycle.
\end{conj} 
P\'osa~(see \cite{posa}) had earlier proposed the conjecture in the case of the square of a Hamilton cycle (that is, when $r=2$). 
Koml\'os, S\'ark\"ozy and Szemer\'edi~\cite{kss} proved Conjecture~\ref{seymour} for graphs whose order is sufficiently larger than $r$. More recently,
Ch\^au, DeBiasio and Kierstead~\cite{cdk} proved P\'osa's conjecture for graphs of order at least $2 \times 10^8$.

In the case when $r+1$ divides $|G|$, a necessary condition for a graph $G$ to contain the $r$th power of a Hamilton cycle is that
$G$ contains a perfect $K_{r+1}$-packing.
Further, notice that the minimum degree condition in Conjecture~\ref{seymour} is the same as the condition in the Hajnal-Szemer\'edi theorem
with respect to perfect $K_{r+1}$-packings. 
Thus an obvious question is whether the condition in Conjecture~\ref{conj1} 
forces a graph to contain the $(r-1)$th power of a Hamilton cycle. Interestingly though, when $r=3$, this is not the case.
\begin{prop}\label{square}
Suppose that $C,n \in \mathbb N$ such that $C \ll n$ and $3$ divides $n$.
Then there exists a graph $G$ whose degree sequence $d_1 \leq \dots \leq d_n$ satisfies
$$d_i \geq \frac{n}{3}+C+i \ \text{ for all } 1\leq i \leq \frac{n}{3} $$
but such that $G$ does not contain the square of a Hamilton cycle.
\end{prop}
\proof Choose $C,K,n \in \mathbb N$ so that $C \ll K \ll n$. 
Let $G$ denote the graph on $n$ vertices consisting of three vertex classes $V_1=\{v\}$, $V_2$ and $V_3$
where $|V_2|=n/3+C+1$ and $|V_3|=2n/3-C-2$ which contains the following edges:
\begin{itemize}
\item All edges from $v$ to $V_2$;
\item All edges between $V_2$ and $V_3$ and all possible edges in $V_3$;
\item There are $K$ vertex-disjoint stars in $V_2$, each of size $\lfloor |V_2|/K \rfloor$,
$\lceil |V_2|/K \rceil$, which cover all of $V_2$ (see Figure~1).
\end{itemize}
\begin{figure}\label{picture}
\begin{center}\footnotesize
\includegraphics[width=0.5\columnwidth]{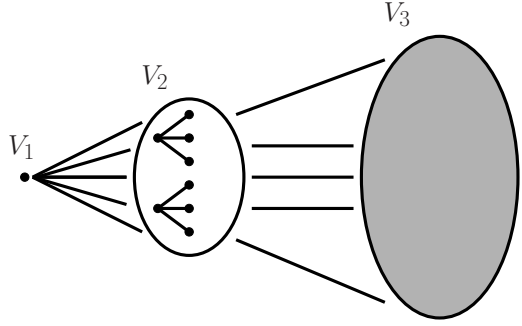}
\caption{The example from Proposition~\ref{square} in the case when $K=2$ and $|V_2|=8$.}
\end{center}
\end{figure}
Let $d_1 \leq \dots \leq d_n$ denote the degree sequence of $G$. There are $n/3+C-K+1 \leq n/3-2C-1 $ vertices in $V_2$
of degree $2n/3-C$. Since $C \ll K \ll n$, the remaining 
$K$ vertices in $V_2$ have degree at least $2n/3-C-2+\lfloor |V_2|/K\rfloor \geq 2n/3 +C+1$.
Since $d_G (v)=n/3+C+1$ and $d_G (x)=n-2$ for all $x \in V_3$, we have that
$d_i \geq \frac{n}{3}+C+i$ for all  $1\leq i \leq \frac{n}{3} $.

A necessary condition for a graph $G$ to contain the square of a Hamilton cycle is that, for every $x \in V(G)$,
$G[N(x)]$ contains a path of length $3$. Note that $N(v)=V_2$ and $G[V_2]$ does not contain a path of length $3$. So $G$
does not contain the square of a Hamilton cycle.
\endproof
Notice that we can set $C=o(\sqrt{n})$ in Proposition~\ref{square}. We finish by raising the following question.
\begin{ques} What can be said about degree sequence conditions which force a graph to contain the $r$th power of a Hamilton
cycle? In particular, can one establish a degree sequence condition that ensures a graph $G$ on $n$ vertices 
contains the $r$th power of a Hamilton cycle and which allows for ``many'' vertices of $G$ to have degree ``much less'' than $rn/(r+1)$?
\end{ques}
\COMMENT{More general/better example?}

\subsection*{Acknowledgements}
We thank the referees for their comments. In particular, we thank one referee for pointing out the work in~\cite{frank}, \cite{dirac} and~\cite{just}.

This research was carried out whilst the third author was visiting the Department of Mathematics 
of the University of Illinois at Urbana-Champaign. 
This author would like to thank the department for the hospitality he received.
We would also like to thank Hal Kierstead for helpful discussions.



\begin{thebibliography}{10}
\bibitem{frank} A. Akiyama and P. Frankl, On the Size of Graphs with Complete-Factors,
\emph{J. Graph Theory}~{\bf 9} (1985), 197--201.
\bibitem{abhp} P. Allen, J. B\"ottcher, J. Hladk\'y and D. Piguet, A density Corr\'adi-Hajnal theorem, \emph{Electronic Notes
in Discrete Mathematics}~{\bf 38} (2011), 31--36.
\bibitem{cdk} P. Ch\^au, L. DeBiasio and H.A. Kierstead, P\'osa's Conjecture for graphs of order at least $2\times 10^8$,
\emph{Random Structures and Algorithms}~{\bf 39} (2011), 507--525.
\bibitem{ch} V. Chv\'atal, On Hamilton's ideals, \emph{J. Combin. Theory~B}~{\bf 12} (1972), 163--168.

\bibitem{corradi} K. Corr\'adi and A. Hajnal, On the maximal number of independent circuits in a graph, \emph{Acta Math. Acad. Sci.
Hungar.}~{\bf 14} (1964), 423--439.
\bibitem{dirac} G.A. Dirac, 
Structural properties and circuits in graphs, in: C.St.-J.A. Nash-Williams, J. Sheehan (Eds.), Proceedings of the 5th British Combinatorial Conference, Congressus Numerantium, No. XV, \emph{Utilitas Math.}, Winnipeg, Man., 1976, pp. 135--140. 
\bibitem{posa} P. Erd\H{o}s, Problem 9, in: M. Fieldler (Ed.), \emph{Theory of Graphs and its Applications}, Czech. Acad. Sci. Publ., Prague, 1964, p. 159.
\bibitem{erdos1} P. Erd\H{o}s, On the graph theorem of Tur\'an, \emph{Mat. Lapok.}~{\bf 21} (1970), 249--251.
\bibitem{hs} A. Hajnal and E.~Szemer\'{e}di, Proof of a conjecture of Erd\H{o}s,
\emph{Combinatorial Theory and its Applications vol. II}~{\bf 4} 
(1970), 601--623.
\bibitem{hell} P. Hell and D.G. Kirkpatrick, On the complexity of general graph
factor problems, \emph{SIAM J. Computing}~{\bf 12} (1983), 601--609.
\bibitem{just} P. Justesen, On independent circuits in finite graphs and a conjecture of Erd\H{o}s and P\'osa,
\emph{Ann. Discrete Math.}~{\bf 41} (1989), 299--306
\bibitem{kier} H.A. Kierstead and A.V. Kostochka, An Ore-type Theorem on Equitable
Coloring, \emph{J. Combin. Theory B}~{\bf{98}}  (2008), 226--234.
\bibitem{knox} F. Knox and A. Treglown, Embedding spanning bipartite graphs of small bandwidth, \emph{Combin. Probab. Comput.}~{\bf 22}
(2013), 71--96.
\bibitem{kss} J. Koml\'os, G.N. S\'ark\"ozy and E. Szemer\'edi, Proof of the Seymour conjecture for large graphs, \emph{Annals of Combinatorics}~{\bf 2} (1998), 43--60.
\bibitem{kuhn} D. K\"{u}hn and D. Osthus, Critical chromatic number and the
complexity of perfect packings in graphs, \emph{17th ACM-SIAM Symposium on 
Discrete Algorithms} (SODA 2006), 851--859.
\bibitem{kuhn2} D. K\"{u}hn and D. Osthus, The minimum degree threshold for perfect
graph packings, \emph{Combinatorica}~{\bf 29} (2009), 65--107.
\bibitem{kotore} D. K\"{u}hn, D. Osthus and A. Treglown, An Ore-type theorem for perfect packings in graphs,
\emph{SIAM J. Disc. Math.}~{\bf 23} (2009), 1335--1355.
\bibitem{seymour} P. Seymour, Problem section, in \emph{Combinatorics: Proceedings of the British Combinatorial
Conference 1973} (T.P. McDonough and V.C. Mavron eds.), 201--202, Cambridge
University Press, 1974.

\end{thebibliography}

\section*{Appendix}
Here we give proofs of Proposition~\ref{propex3} and Lemma~\ref{case3prop}. The following fact will be used in both of these proofs.
\begin{fact}\label{difffact}
Fix $n,r \in \mathbb N$ such that $r\geq 3$ and $r$ divides $n\geq 2r$. Define
$$h(x):=x+\frac{(n-x-1)^2}{2(r-2)}-\frac{1}{2}(n-x-1).$$
Then $h(x)$ is a decreasing function for $ x \in [0, n/(r-1)]$.
Moreover, if $n\geq 3r$ then $h(x)$ is a decreasing function for $ x \in [0, (n+r)/(r-1)]$.
\end{fact}
\proof
Notice that
$$h'(x)=\frac{3}{2}-\frac{(n-x-1)}{r-2} =\frac{x}{r-2}+\frac{1-n}{r-2}+\frac{3}{2}.$$
So for $x \leq n/(r-1)$,
\begin{align*}
h'(x)  \leq \frac{n}{(r-1)(r-2)} +\frac{1-n}{r-2}+\frac{3}{2}= -\frac{n}{r-1}+\frac{1}{r-2}+\frac{3}{2}.
\end{align*}
Note that $3(r-1)/2+(r-1)/(r-2) < n$ since $n\geq 2r$ and $r \geq 3$. Thus,
$$h'(x) \leq -\frac{n}{r-1}+\frac{1}{r-2}+\frac{3}{2} <0.$$
If $x \leq (n+r)/(r-1)$ then
$$h'(x) \leq \frac{n+r}{(r-1)(r-2)} +\frac{1-n}{r-2} +\frac{3}{2}=-\frac{n}{r-1}+\frac{1}{r-2}+ \frac{r}{(r-1)(r-2)}+\frac{3}{2}.$$
If $n \geq 3r$ then $n >3r/2+4$. So $n> 3(r-1)/2+(2r-1)/(r-2)$. Thus, 
$$h'(x) \leq -\frac{n}{r-1}+\frac{1}{r-2}+\frac{r}{(r-1)(r-2)}+\frac{3}{2} <0,$$
as desired.
\endproof
\medskip
\noindent
{\bf Proof of Proposition~\ref{propex3}.}
We need to show that, for all $D \in \mathbb N$ such that $n/r \leq D \leq n/(r-1)$,
$$\frac{n^2}{2r^2}+\frac{n}{2r} =\binom{n/r+1}{2}   \leq D+e(\overline{T}(n-D-1, r-2)) .$$
Since $D \leq n/(r-1)$, Facts~\ref{turanbound} and~\ref{difffact} imply that
\begin{align*}
D+e(\overline{T}(n-D-1, r-2)) & \geq D+ \frac{(n-D-1)^2}{2(r-2)}-\frac{(n-D-1)}{2}  \\ &
\geq \frac{n}{r-1}+\frac{1}{2(r-2)}\left[ \frac{(r-2)}{r-1}n- 1 \right] ^2-\frac{1}{2}\left[\frac{(r-2)}{r-1}n- 1\right] \\ &
\geq \frac{(r-2)}{2(r-1)^2} n^2 -\frac{(r-2)}{2(r-1)}n.
\end{align*}
Thus, it suffices to show that
\begin{align}\label{target}
\frac{(r-2)}{2(r-1)^2} n -\frac{r-2}{2(r-1)} \geq \frac{n}{2r^2}+\frac{1}{2r}.
\end{align}
Notice that 
\begin{align}\label{lemhelp}
\frac{r-2}{2(r-1)^2} -\frac{1}{2r^2} =\frac{(r-2)r^2-(r-1)^2}{2r^2(r-1)^2}=\frac{r^3-3r^2+2r-1}{2r^2(r-1)^2}
\end{align}
and 
$$\frac{r-2}{2(r-1)}+\frac{1}{2r}=\frac{r^2-r-1}{2r(r-1)}.$$
Since $n \geq 2r$, (\ref{target}) implies that it suffices to show that
\begin{align}\label{target2}
\frac{r^3-3r^2+2r-1}{r(r-1)^2} -\frac{r^2-r-1}{2r(r-1)}\geq 0.
\end{align}
Note that $r^3\geq 4r^2-4r+3$ as $r\geq 3$. Thus,
$2(r^3-3r^2+2r-1)\geq (r^2-r-1)(r-1)$. So indeed (\ref{target2}) is satisfied,
as desired.
\endproof
\medskip
\noindent
{\bf Proof of Lemma~\ref{case3prop}.} We need to show that, for all $D \in \mathbb N$ such that $n/r \leq D < (n+r)/(r-1)$,
$$\frac{n^2}{2r^2}+\frac{n}{2r} =\binom{n/r+1}{2}   \leq D+e(\overline{T}(n-D-1, r-2)) .$$
Since $D<(n+r)/(r-1)$ we have that $D\leq n/(r-1)+1$. So Facts~\ref{turanbound} and~\ref{difffact} imply that
\begin{align*}
D+e(\overline{T}(n-D-1, r-2)) & \geq D+ \frac{(n-D-1)^2}{2(r-2)}-\frac{(n-D-1)}{2}  \\ &
\geq \frac{n}{r-1}+1+ \frac{1}{2(r-2)}\left[ \frac{(r-2)}{r-1}n- 2 \right] ^2-\frac{1}{2}\left[\frac{(r-2)}{r-1}n- 2\right] \\ &
\geq \frac{(r-2)}{2(r-1)^2} n^2 -\frac{(r-2)}{2(r-1)}n- \frac{n}{r-1}.
\end{align*}
Thus, it suffices to show that 
\begin{align}\label{targety}
\frac{(r-2)}{2(r-1)^2} n -\frac{(r-2)}{2(r-1)}-\frac{1}{r-1} \geq \frac{n}{2r^2}+\frac{1}{2r}.
\end{align}
Notice that $$\frac{r-2}{2(r-1)}+\frac{1}{r-1}+\frac{1}{2r}=\frac{r^2+r-1}{2r(r-1)}.$$
Since $n\geq 3r$, (\ref{lemhelp}) and (\ref{targety}) imply that it suffices to show that
\begin{align}\label{targety2}
\frac{3(r^3-3r^2+2r-1)}{2r(r-1)^2} -\frac{r^2+r-1}{2r(r-1)}\geq 0.
\end{align}
Note that $2r^3-9r^2+8r-4 \geq 0$ as $r\geq 4$. Thus, $3(r^3-3r^2+2r-1)\geq (r^2+r-1)(r-1)$.
So indeed (\ref{targety2}) is satisfied, as desired.
\endproof

\end{document}